\documentclass[12pt]{article}
\usepackage{amsmath,amsxtra}
\usepackage{amssymb}
\usepackage{txfonts}

\usepackage{diagrams}

\newtheorem{theorem}{Theorem}[section]
\newtheorem{deff}[theorem]{Definition}

\newtheorem{proposition}[theorem]{Proposition}

\newtheorem{lemma}[theorem]{Lemma}

\numberwithin{equation}{section}

\newcommand{\npa}{\addtocounter{theorem}{1} \noindent
{\bf \arabic{section}.\arabic{theorem}}\;\;}

\newcommand{\cU}{{\cal U}}
\newcommand{\cG}{{\cal G}}

\newcommand{\proof}{\vskip2mm \noindent {\bf Proof.}~}

\newcommand{\mto}{\mapsto}
\newcommand{\na}{\nabla}

\newcommand\beqa {\begin{eqnarray}}
\newcommand\eeqa {\end{eqnarray}}
\newcommand\bqa {\begin{eqnarray}}
\newcommand\eqa {\end{eqnarray}}
\newcommand{\beq}{\begin{eqnarray}}
\newcommand{\beqn}{\begin{eqnarray}\nonumber}
\newcommand{\eeq}{\end{eqnarray}}
\newcommand{\be}{\begin{array}}
\newcommand{\ee}{\end{array}}

\newcommand{\e}{\epsilon}

\newcommand\bea {\begin{eqnarray}}
\newcommand\eea {\end{eqnarray}}

 \newcommand{\pt}{\partial}

 \newcommand{\cS}{{\cal S}}

 \newcommand{\cP}{{\cal P}}
 
 \newcommand{\cA}{{\cal A}}
 \newcommand{\cD}{{\cal D}}

 \newcommand{\cQ}{{\cal Q}}
 
 \newcommand{\cM}{{\cal M}}
 \newcommand{\cL}{{\cal L}}
 \newcommand{\cN}{{\cal N}}
 \newcommand{\cF}{{\cal F}}
 
 \newcommand{\cI}{{\cal I}}

\newcommand{\R}{{\mathbb R}}
\newcommand{\Z}{{\mathbb Z}}

 \newcommand{\g}{{\mathfrak g}}

 \newcommand{\h}{{\mathfrak h}}

\newcommand{\md}{\mathrm{d}}

 \def\S{{\Sigma}}
\def\2{{\textstyle\frac{1}{2}}}
\def\ba{\begin{eqnarray}}
\def\ea{\end{eqnarray}}

\newcommand{\TS}[1]{{\bf\boldmath \framebox{TS} :#1}}

   \def\CM{{\cal M}}

\makeatletter
\def\bard{\protect\@bard}
\def\@bard{%
\relax \bgroup
\def\@tempa{\hbox{\raise.73\ht0
\hbox to0pt{\kern.4\wd0\vrule width.7\wd0
height.1pt depth.1pt\hss}\box0}}%
\mathchoice{\setbox0\hbox{$\displaystyle\mathrm{d}$}\@tempa}%
{\setbox0\hbox{$\textstyle \mathrm{d}$}\@tempa}%
{\setbox0\hbox{$\scriptstyle \mathrm{d}$}\@tempa}%
{\setbox0\hbox{$\scriptscriptstyle \mathrm{d}$}\@tempa}%
\egroup } \makeatother

\makeatletter
\def\barp{\protect\@barp}
\def\@barp{%
\relax \bgroup
\def\@tempa{\hbox{\raise.73\ht0
\hbox to0pt{\kern.4\wd0\vrule width.7\wd0
height.1pt depth.1pt\hss}\box0}}%
\mathchoice{\setbox0\hbox{$\displaystyle\partial$}\@tempa}%
{\setbox0\hbox{$\textstyle \partial$}\@tempa}%
{\setbox0\hbox{$\scriptstyle \partial$}\@tempa}%
{\setbox0\hbox{$\scriptscriptstyle \partial$}\@tempa}%
\egroup } \makeatother

\makeatletter
\def\barb{\protect\@barb}
\def\@barb{%
\relax \bgroup
\def\@tempa{\hbox{\raise.73\ht0
\hbox to0pt{\kern-.1\wd0\vrule width.7\wd0
height.1pt depth.0pt\hss}\box0}}%
\mathchoice{\setbox0\hbox{$\displaystyle\mathrm{b}$}\@tempa}%
{\setbox0\hbox{$\textstyle \mathrm{b}$}\@tempa}%
{\setbox0\hbox{$\scriptstyle \mathrm{b}$}\@tempa}%
{\setbox0\hbox{$\scriptscriptstyle \mathrm{b}$}\@tempa}%
\egroup } \makeatother

   \def\a{\alpha}
   \def\b{\beta}

   \def\de{\delta}

   \def\e{\epsilon}
   \def\i{\imath}

\def\Gr{\mathrm{G}} 

\def\G1{\Gr_1}

\def\dot{{\! \bullet}}

\def\weco{\stackrel{\wedge}{,}}

\def\gf{\varphi}
\def\map{m}

\begin{document}

\begin{center}
{ \bf \Large Characteristic classes associated to Q-bundles}
\end{center}


\begin{center}\sl\large
Alexei Kotov \footnote{e-mail address: {\tt Alexei.Kotov @
uni.lu}}, Thomas Strobl \footnote{e-mail address: {\tt Strobl @
math.univ-lyon1.fr}}
\end{center}

   \begin{center}{\small\sl $^1$ Laboratoire de  Mathematiques,
    Universite du Luxembourg \\
  \vspace{0.2em}
$^2$ Institut Camille Jordan, Universite Claude Bernard Lyon 1}
\end{center}

\vskip 3mm

{\small \noindent {\bf Abstract:} A Q-manifold is a graded manifold
endowed with a vector field of degree one squaring to zero. We
consider the notion of a Q-bundle, that is, a fiber bundle in the category
of Q-manifolds. To each homotopy class of ``gauge fields'' (sections in
the category of graded manifolds) and each cohomology class of a
certain subcomplex of forms on the fiber we associate a cohomology
class on the base. Any principal bundle yielding canonically a
Q-bundle, this construction generalizes Chern-Weil classes.
Novel examples include cohomology classes that are locally the de
Rham differential of the integrands of topological sigma models
obtained by the AKSZ-formalism in arbitrary dimensions. For
Hamiltonian Poisson fibrations one obtains a characteristic
3-class in this manner. We also relate to equivariant cohomology
and Lecomte's characteristic classes of exact sequences of Lie
algebras.



\vskip 4mm


\vspace{0.2em}

}

\section{Introduction}

\npa The notion of a Q-manifold provides a general framework for
studying gauge theories within the  Batalin-Vilkovisky formalism
\cite{Schwarz_semiclassical}.  It is particularly useful in the context of
topological sigma models \cite{AKSZ}. A Q-manifold, also known as a
differential-graded (dg) manifold, is a graded manifold $\cM$ endowed
with a degree one vector field $Q$ which satisfies the equation
$[Q,Q]\equiv 2 Q^2=0$. Hereafter we suppose that the algebra of functions on
$\cM$ is non-negatively graded, unless the contrary is
stated. We say that $\cM$ is a Qp-manifold, if the algebra of
functions is locally generated in degree up to p.

\vskip 3mm

\npa Let us enumerate some basic examples of Q-manifolds
 appearing in the literature.

\vskip 2mm

\noindent  (1) A Lie algebra $\g$ considered as a purely odd
manifold of degree one, denoted as $\g [1]$. The algebra of
functions is naturally isomorphic to $\Lambda^\dot \g^*$, the cochain
complex of $\g$, and the Q-field is nothing but the
Chevalley-Eilenberg differential.

\vskip 2mm

\noindent (2) A Lie algebroid $E$ with the degree of fibers
shifted by one, denoted as $E [1]$. The algebra of functions is
identified with $\Gamma (\Lambda^\dot E^*)$ and the Q-field is the
canonical differential. Moreover, every Q1-manifold is necessarily
of the form $E[1]$ for a certain Lie algebroid \cite{Vaintrob}.
In general, a homological vector field of degree one on an arbitrary
graded vector bundle  determines an $L_\infty -$algebroid
structure by use of multi-derived brackets, cf.,~e.g.,~\cite{Voronov_derived}.

\vskip 2mm

 \noindent (3) A PQ-manifold: This is a graded manifold $\cS$ supplied
 with a symplectic form $\omega$ of degree p and a function $\cQ$ of
 degree p+1, which obeys the equation of self-commutativity with
 respect to the non-degenerate Poisson bracket determined by
 $\omega$. The Q-field is the hamiltonian vector field of $\cQ$.\footnote{For
 $p>1$ a PQ-manifold is \emph{equivalent} to a Q-manifold with compatible
 degree p symplectic form \cite{Roytenberg0203110}.}

\vskip 3mm

\npa
A morphism of Q-manifolds (Q-morphism) is a degree preserving
map $\phi$, the pull-back of which commutes with the corresponding
homological vector fields, considered as super derivations of
functions, i.e. the following chain property holds: $ Q_1\phi^*
=\phi^* Q_2 $. A morphism of Q1-manifolds is nothing but the
morphism of the corresponding Lie algebroids (\cite{Vaintrob}; cf.~\cite{BKS} for
a proof of equivalence with the original definition of Lie algebroid morphisms
given in \cite{Mac}). Given a smooth map of two manifolds $\map
\colon  M\to N$, its push-forward defines a Q-morphism $\map_* \colon  T[1]M\to
T[1]N$ of the tangent bundles, where a Q-structure on the odd tangent
bundle of a manifold is determined by the de Rham operator
regarded as a homological vector field by use of the
identification $C^\infty (T[1]M)\simeq \Omega^\dot (M)$.

\vskip 2mm\noindent Apparently, a composition of two Q-morphisms
is again a Q-morphism, so there is a well-defined category of
Q-manifolds. A Q-bundle, a \emph{fiber} bundle in this category, is a
surjective morphism of the total space to the base of the bundle,
satisfying an additional requirement of local triviality: a bundle
is built from direct products of local charts on the base and a
fixed fiber glued by a transition cocycle of ``gauge
transformations''---as will be  detailed further in section 2 below.
By a ``gauge field'' $\gf$ in a Q-bundle $\pi
\colon  \cM\to\cM_1$ we mean a section of the underlying bundle of graded
manifolds. In general we do not assume that a section is a
Q-morphism! (The existence of such a section imposes a certain
constraint on the bundle).

Some examples of these constructions are the following ones, with
the third one providing the relation to ordinary gauge theories:

\vskip 2mm\noindent (1) The product of two Q-manifolds is again a
Q-manifold and the projection to each factor produces a (trivial)
Q-bundle structure.

\vskip 2mm\noindent (2) A fiber bundle $p\colon  M\to M_1$ determines a
natural ``non-linear'' example of a Q-bundle by use of the the
push-forward map: $p_*\colon  T[1]M\to T[1]M_1$. Apparently, the
push-forward of any section of $p$ is a section of $p_*$, which
is, indeed, a morphism of the corresponding Q-manifolds.

\vskip 2mm\noindent (3) Given a principal G-bundle $p\colon P\to M$, we
construct a Q-bundle in the following way: As total space we take the
quotient of $T[1]P /G$, where the group action of $G$ on $P$ is lifted
in the canonical way and the quotient by $G$ can be considered as a
bundle over $M$. Using the push-forward $p_*$ of $p$, on the other
hand, we obtain a (degree-preserving) map to $T[1]M$, the base of the
Q-bundle. Both spaces are canonically equipped with the de Rham
differential (in the first case restricted to $G$-invariant
differential forms on $P$).  This construction is known as the Atiyah
algebroid of $P$, which is a particular Lie algebroid ($T[1]P$ is
obviously a degree one $Q$-manifold). A connection in $P$ provides a
lift of tangent vectors on $M$ to tangent vectors of $P$; by its
equivariance w.r.t.~the $G$-action this corresponds precisely to a
bundle map $\gf \colon TM \to TP/G$, i.e.~a section of the bundle $p_*
\colon T[1]P /G \to T[1]M$. As we will see in detail in section 2 below,
the connection is flat, iff $\gf$ is a Q-morphism.

\vskip 2mm\noindent (4) A transitive Lie algebroid $E\to M$, in
generalization of an Atiyah algebroid: By definition this is a Lie
algebroid with surjective anchor $\rho$, yielding the short exact
sequence \beq 0\to\underline{\g} \to E \stackrel{\rho}{\to} TM\to
0 \, , \label{seq} \eeq where $\underline{\g}$ is a bundle of Lie
algebras defined by the kernel of $\rho$. We restrict to the case
that any of the fibers $\ker \rho$ is isomorphic to a single Lie
algebra $\g$.\footnote{Note that even under this assumption not
every transitive Lie algebroid comes from a principal bundle. Only
if this Lie algeroid can be integrated to a Lie groupoid, this is
the case.}.  This then yields a Q-bundle $\rho \colon E[1] \to
T[1]M$ with typical fiber $\g[1]$. Note that since $\rho$ is a
morphism of Lie algebroids, the projection is a Q-morphism. In
this particular case, a gauge field $\gf \colon T[1]M\to E[1]$ is
a splitting of the exact sequence of Lie algebroids (\ref{seq}).
It is sometimes also called a ``connection'' of the transitive Lie
algebroid $E$, in generalization of the previous example, and
called ``flat'' in a situation when $\gf$ is a Q-morphism.

\vskip 2mm\noindent (5) An exact sequence of Lie algebras---cf.~example (1) of paragraph
1.2---is a Q-(fiber)-bundle (as defined above) only in the case when
it is isomorphic to a direct sum of Lie algebras. We will address
this situation at the end of the paper.

\vskip 2mm\noindent (6) More generally than examples (2) - (5),
one can consider an exact sequence of Lie algebroids, covering an
ordinary fiber bundle; if the total Lie algebroid splits locally
into a direct product of fiber and base Lie algebroids, it fits
into the definition of Q1-bundles above. We will study particular
examples of this, where the fibers are some given PQ-manifold and
the base a tangent Lie algebroid, considering applications in
section 4 below.

\vskip 3mm\npa For an arbitrary degree preserving map of
Q-manifolds $\gf \colon \cM_1\to\cM_2$ the difference $F:= Q_1
\gf^* -\gf^* Q_2$, which we call the ``field strength'' if $\gf$
is a (coarse-grained) section of a Q-bundle $\cM_2 \to \cM_1$, is
non-vanishing in general. It is a degree one derivation of
functions on the target taking values in functions on the source
\beq F\colon C^\infty (\cM_2)\to C^\infty (\cM_1)\;,\eeq for which
the following Leibnitz-type property holds: \beq \label{Leib} F
(gh)=F(g)\gf^*(h)+(-1)^{\deg(g)}\gf^*(g)F(h)\;,\hspace{3mm}\forall
g,h\in C^\infty (\cM_2)\;.\eeq Therefore $F$ can be identified
with a degree one section of the pull-back bundle $\gf^* (T
\cM_2)$ over $\cM_1$ or, as is equivalent, with a degree
preserving map $f \colon \cM_1\rightarrow T [1]\cM_2$ covering
$\gf$, as will be further detailed in section 3 below. The graded
manifold $T [1]\cM$, where $\cM$ is a Q-manifold, is a double
Q-manifold (or a double Q-algebroid,
cf.~\cite{Mackenzie_double,Voronov_Qman}), i.e.~it admits a pair
of anti-commuting homological vector fields. It will turn out,
cf.~Proposition \ref{f-chain_map} below, that $f$ is a Q-morphism
if $T [1]\cM_2$ is endowed with the sum of two canonical
Q-structures as differential. Note that in contrast to the
Leibnitz property (\ref{Leib}) of $F$, the map $f^*\colon
C^\infty(T[1]\cM_2) \to C^\infty(\cM_1)$ defines a morphism of
algebras, and thus, being a chain map, also a map in cohomologies.

Using the Bernstein-Leites sign convention, functions on $T
[1]\cM$ can be identified with differential forms on $\cM$. Given
a trivialization of a Q-bundle over some open cover of the base,
we can always identify a section with a family of degree
preserving maps from the local charts to the fiber, which are
related by the transition transformations on double overlaps. By
means of the pullback of the above map f, we obtain a family of
chain maps acting from the total complex of differential forms on
the fiber to the complex of functions on the open charts, which
are different in general on double overlaps. However, applying the
collection of chain maps to a differential form on the fiber,
which is invariant with respect to the gluing transition functions
(or, equivalently, the ``gauge transformations'')---we will call
such forms \emph{basic}---, we obtain a well-defined cocycle on
the whole base.  This thus provides a map from the cohomology of
basic forms on the fiber to the cohomology of forms on the base,
 cf.~Theorem \ref{construction_of_charclasses} below. The construction will be
seen to generalize the Chern-Weil formalism of characteristic classes
associated to principal bundles.



\vskip 3mm\npa Some brief remark on our notation and nomenclature:
As we recalled above, a connection in a principal bundle $P$ is in
bijection to a splitting in (\ref{seq}) of the corresponding Atiyah
algebroid $E \to M$, or, what is the same, a (globally well-defined)
section $\gf$ (in the coarse-grained sense) of a likewise Q-bundle $\rho
\colon E[1] \to T[1]M$. The map $f$ above can be seen to
generalize the curvature of that connection on $P$, moreover. Since,
on the other hand, on a Q-bundle $\pi\colon\CM \to \CM_1$ one can also
discuss (super-) connections and curvatures, we refrained from calling
such maps $\gf$ and $f$ as (generalized) connections and curvatures,
respectively. Instead we thus prefer a physics oriented nomenclature
in this context, calling $\gf$s ``gauge fields'' and $f$s ``field
strengths''; in the context of $P$, they locally are represented by
(Lie algebra valued) 1-forms and 2-forms $A$ and $F\sim F_A$,
respectively.
Note that in a
more general situation $\gf$ may correspond to a collection of
differential forms of different degrees, or, when $\CM_1$ is not the
shifted tangent bundle of some manifold, even not to differential forms
at all.

 ``Gauge transformations'' or ``gauge symmetries'' will
turn out to be related also to vertical automorphisms of the Q-bundle
(vertical automorphisms of a principal bundle give rise to anchor
preserving automorphisms of the Atiyah algebroid), but in general it
will be useful to consider only a subset of the latter ones for gluing
transformations.

\vskip 3mm\npa The paper is organized as follows. In section 2 we
discuss the appropriate notion of gauge transformations, generalizing
\cite{BKS}, for a Q-bundle and describe their action on the space of
gauge fields.

\vskip 2mm\noindent In section 3 we prove the chain property of
the map $f$, defined above, and show that the Weil algebra model
of characteristic classes is a particular case of our
construction. Here we slightly adapt the notion of a basic form,
looking at the action of gauge symmetries on $f$, and complete the
construction of characteristic classes associated to a section of
a Q-bundle. We also prove the obtained cohomological classes are
homotopical invariants of such sections.

\vskip 2mm\noindent In section 4 we explain a possible
construction of characteristic classes whose cocycles turn out to
be locally represented by integrands of the (classical part of)
topological AKSZ-type sigma models. In the case of Hamiltonian
Poisson fibrations, the typical fiber being a Poisson manifold,
one obtains a 3-class in de Rham cohomology on the base manifold,
which, locally, agrees with the original construction of (the
integrand of) the Poisson sigma model in \cite{Schaller-Strobl}.
In section 4 we also address the gauging of Wess-Zumino terms, its
relation to equivariant cohomology, and the characteristic classes
of \cite{Lecomte} within the present framework.

\section{Q-bundles and gauge symmetries}
\label{sectionQbundles}

\npa In order to explain the notion of gauge fields and gauge
symmetries, we start with a simple example of $\g -$valued 1-forms
on a smooth manifold $M$, where $\g$ is a Lie algebra. Given
$A\in\Omega^1 (M,\g )$, interpreted as a connection in a trivial
bundle $M \times G$, Lie($G$)$=\g$, we look at its curvature: \beq F_A :=\md
A +\frac{1}{2}[A,A]\;.\eeq The group of $G-$valued functions on
$M$ is acting on connections by:
$A^g =g^{-1}\md g + \mathrm{Ad}_{g^{-1}}
(A)$, where $g^{-1}\md g$ is the pull-back by $g$ of the (left)
Maurer-Cartan form on the Lie group and $\mathrm{Ad}$ is the
adjoint action. 
The above transformations correspond to vertical automorphisms of the
trivial bundle and are called gauge transformations in the physics
literature. Their infinitesimal version is goverend by a $\g-$valued
function $\e$: \beq\label{gauge_linear} \de_\e A:=\frac{\md}{\md t}
A^{\exp(t\e)}\mid_{t=0}=\md\e +[A,\e]\;.\eeq
The condition of flatness,
$F_A=0$, can be also regarded as Maurer-Cartan equation for $A$.

\vskip 2mm\npa Let us adapt this example to the language of
dg or Q-manifolds. As we already know, a Lie
algebra can be treated as a Q-manifold $\g [1]$, such that the
algebra of functions becomes isomorphic to $\Lambda(\g^*)$ with
the Q-field given by  the Chevalley-Eilenberg differential:
\beq\label{Chevalley-Eilenberg}\md_\g
(\a)(\eta, \eta ')=-\a ([\eta, \eta'])\eeq where $\a\in\g^*$
and $\eta ,\eta'\in\g$. The product of $T[1]M$ and $\g [1]$ is
again a Q-manifold, the Q-structure of which is given by the sum
of de Rham and Chevalley-Eilenberg derivations extended to the
product in the standard way. A $\g-$valued 1-form on $M$ can be
thought of as a degree preserving map $\gf \colon  T[1]M \to \g [1]$
and its graph as a section of the bundle
\beq\label{Q-bundle_trivial} T[1]M\times \g [1]\to T[1]M\;.\eeq The
pull-back of $\gf$ is acting as follows: for each
$\omega\in\Omega (M)$, $\a\in\Lambda^p (\g^*)$ one has \beqn
\gf^* (\a\otimes\omega)=\a
(\,\underbrace{A\weco \,\ldots\,\weco A}_{p\,\mathrm{times}}\,)\wedge\omega\;.\eeq
Any $\g-$valued function $\e$, acting by the contraction $\iota_\e$ on
$\Omega (M)\otimes \Lambda (\g^*)$, can be considered as a
super-derivation of degree -1, which super-commutes with $\Omega
(M)$. The last property implies that it can be identified with a
vertical vector field on the total space of
(\ref{Q-bundle_trivial}).
\begin{proposition} \label{trivbundle} The following identity holds
for each $\omega\in\Omega (M)$ and $\a\in\Lambda^p (\g^*)$:
\beq\label{Linear_curvature} \left(\md \,\gf^* -\gf^* \,(\md +
\md_\g)\right)\, (\a\otimes\omega ) &=& \sum_k (-1)^{k+1} \a
(\,A \weco \,\ldots \weco
\overbrace{F_A}^{k} \weco \ldots \weco A\,)\wedge\omega,\\
\label{gauge_vs_derivative}\gf^* L_\e \, (\a\otimes\omega) &=&
\sum_k \a (\,A,\weco\,\ldots \weco \overbrace{\de_\e A}^{k} \weco \ldots \weco
A\,)\wedge\omega\;,\eeq where $L_\e$ is the Lie derivative along
$\e$, defined as the super-commutator $[Q, \iota_\e ]$ and $\de_\e A$ is given by formula
(\ref{gauge_linear}) above.
\end{proposition}

\proof Straightforward calculations. $\square$

\vskip 2mm\noindent Note that, instead of thinking of an
infinitesimal gauge transformation as a flow on the space of
connections, we define a vector field on the total space of
(\ref{Q-bundle_trivial}), the action of which on the space of
connections, regarded as sections of (\ref{Q-bundle_trivial}), can
be naturally induced. One may wonder why, though $A$ is extended
as a morphism of graded manifolds, its infinitesimal variation is
extended (by the Leibnitz rule) as a derivation. Indeed, it is a
general fact, adapted to the graded case, that the space of
infinitesimal variations (the tangent space) of a smooth map $\psi
\colon M\to N$ can be identified with the space of sections of the
pullback bundle $\psi^* (TN)$ or, equivalently, with the space of
derivations $\de \colon C^\infty (N)\to C^\infty (M)$ covering $\psi$:
\beqn\de (hh')=\de (h)\psi^* (h')+ (-1)^{\deg (\de)\deg (h)}\psi^*
(h)\de (h')\;\eeq for any $h,h' \in C^\infty (N)$.  As it is clear
from (\ref{Linear_curvature}), the curvature $F_A$ is the only
obstruction for $A$ to be a Q-morphism; $A$ gives a flat connection,
$F_A=0$, iff the corresponding section $\gf \colon T[1]M \to T[1]M
\times \g[1]$ is a Q-morphism.

\vskip 2mm\npa It is evident how to generalize the picture
described above for connections in a non-trivial vector bundle or its
associated principal bundle $P$: For this purpose we return to example
(3) in paragraph 1.3 above, replacing the trivial Q-bundle
(\ref{Q-bundle_trivial}) by its Atiyah algebroid $E$, i.e.~by
\beq\label{Q-bundle_linear}\rho\colon  E[1]\to T[1]M\,,\eeq where $E=TP/G$ and
$\rho = p_*$ denotes the anchor map. Then a connection becomes a
section of (\ref{Q-bundle_linear}). An infinitesimal gauge
transformation a Lie derivative with respect to some section of $E$
lying in the kernel of $\rho$. The space of such sections is in
one-to-one correspondence with the space of vertical vector fields of
degree minus one on the Q-bundle $E[1] \to T[1]M$; the correspondence
is given by the contraction, $\iota_\epsilon$ for any $\epsilon \in
\Gamma(E)$ can be regarded as a vector field on $E[1]$. The gauge
transformations are thus generated by $L_\e\equiv [Q,\iota_\epsilon]$,
where $\e \in \Gamma(\ker \rho)$ or, equivalently, $\rho_* \iota_\e
=0$. Elements $\alpha \otimes \omega\in \Lambda^p (\g^*) \otimes
\Omega (M)$, used in Proposition \ref{trivbundle} above now generalize
merely to functions on $E[1]$.

\vskip 2mm\npa Now we can describe a general Q-bundle, its gauge
transformations and fields. The Lie super-algebra of vector fields
on $\cM$, denoted as $\cD (\cM)=\oplus_k \cD^k (\cM)$, is a
differential graded Lie algebra, the differential of which is
given by the adjoint action of $Q$: $ad_Q (X):=[Q,X]$ for $X\in
\cD (\cM)$. Indeed, since $Q$ satisfies the ``master equation''
$[Q,Q]=0$, its double commutator with any vector field vanishes:
by the super Jacobi identity one has \beqn ad_Q^2
(X)=[Q,[Q,X]]=\frac{1}{2}\left[[Q,Q],X\right]\equiv 0\;.\eeq We
call a vector field $X\in \cD^0 (\cM)$ commuting with $Q$ an {\em
infinitesimal symmetry}, or simply symmetry, of a Q-manifold and a
degree zero vector field, which is a commutator of $Q$ with some
other vector field, an {\em inner derivation} or inner
(infinitesimal) symmetry. In the case of a Lie algebra $\g[1]$
this agrees with the usual nomenclature, whereas in the case of
$T[1]M$ e.g.~all symmetries are inner in this sense and correspond
to vector fields, i.e.~(infinitesimal) diffeomorphisms. By
definition, infinitesimal symmetries and inner derivations are
cocycles and coboundaries in $(\cD (\cM),ad_Q)$, respectively. The
inclusion ``coboundaries $\subset$ cocycles'' implies that all
inner derivations are infinitesimal symmetries of a Q-manifold.
The following identity follows from the super Jacobi identity and
the nilpotency of $ad_Q$: \beq\label{commutator_Lie} [ad_Q (X),
ad_Q (X')]=ad_Q \left([X,X']_Q\right)\;,\eeq where
$[X,X']_Q:=(-1)^{\deg X+1}[ad_Q(X), X']\equiv[[X,Q],X']$ is known
in mathematics as the derived bracket \cite{Kosmann}. The space of
vector fields supplied with the derived bracket is an example of
a Loday algebra. Note that the
derived bracket between two vector fields is not (super)
skew-symmetric, unless the vector fields are super-commuting. Obviously,
vector fields of degree minus one are closed with respect to the
derived bracket. Moreover, as is clear from
(\ref{commutator_Lie}), but also can be verified directly, its
image by $ad_Q$ is a (super) \emph{Lie} subalgebra in $\cD^0
(\cM)$.


\begin{deff} A subgroup $H$ of degree preserving maps is called a
subgroup of automorphisms (respectively, inner automorphisms), if
its Lie algebra consists of infinitesimal symmetries
(respectively, inner infinitesimal symmetries).
\end{deff}

\noindent Before writing the general definition, let us examine
once more the trivial example, which is a brick underlying a global
design. Suppose $\cM=\cN\times\cF$ is a product of two Q-manifolds
$\cN$ and $\cF$ and $\pi\colon \cM\to\cN$ is a bundle given by the
projection to the first factor. It is obvious that the space of
vertical vector fields can be identified with sections of the
pull-back of $T\cF$ w.r.t.~the second projection.

\begin{proposition}\label{derived_on_total} Let $\cG$ be a graded
Lie subalgebra of vector fields on $\cF$, closed under the derived
bracket. Then the space of functions on $\cN$ taking values in
$\cG$ is a Lie subalgebra of vertical vector fields closed under
the derived bracket on the total space $\cN \times \cF$.
\end{proposition}
\proof Let us take an arbitrary $Y\in C^\infty (\cN, \cG)$, which
can always be written as a linear combination $Y=\sum_j b^j Y_j$,
where $b^j$ are functions on the base and $Y_j\in\cG$. Then the
inner derivative generated by $Y$ on the total space is
\beq\label{ad_on_total_space} ad_Q (Y)=[Q_1 +Q_2, Y]=\sum_j Q_1
(b^j )Y_j +(-1)^{\deg(b^j)}b^j ad_{Q_2}(Y_j)\;.\eeq For any
$X=\sum_i a^i X_i$, its derived bracket with $Y$ is again a
function taking values in $\cG$: \beqn [X,Y]_Q = \sum_{i,j}
(-1)^{\deg(X_i)+\deg(a^i)+1}\left( [Q_1 (a^i )X_i
+(-1)^{\deg(a^i)}a^i ad_{Q_2}(X_i), b^j Y_j]\right) =\\\nonumber
\sum_{i,j}\left( (-1)^{\deg(a^i)+\deg(X_i)(\deg(b^j)+1)+1} Q_1
(a^i )b^j [X_i, Y_j]+ (-1)^{(\deg(X_i)+1)\deg(b^j)}a^i b^j
[X_i,Y_j]_Q\right) \;,\eeq which implies that $\cG$ is a Loday
algebra with respect to the total Q-structure. $\square$

\vskip 2mm\noindent Let us use the notation for the following Lie
algebra of vector fields on the total space:
 \beq\label{def_of_G-prime} \cG' :=ad_Q
\left(C^\infty (\cN, \cG)\right)\cap \cD^0 (\cN \times \cF)\;.\eeq
Is is not
a surprise for us that $\cG'$ consists of vertical vector fields.
Indeed, suppose we are given $X\in\cG'$, then there exists some
element $\e\in\cG$ such that $X=[Q,\e]$. Both of two vector fields
in the commutator are $\pi -$projectable, since $\pi_* (Q)=Q_1$
and $\pi_* (\e)=0$, thus $\pi_* (X)=0$. It well-known that
exponentiating a vertical vector field (at least locally), we
obtain a fiber-wisely acting automorphism, i.e. an automorphism
$\Psi$ satisfying $\pi\circ\Psi =\pi$. Apparently, the set of
fiber-wisely acting automorphisms is a subgroup of all
automorphisms of a bundle and a composition of $\Psi$ with any
section of $\pi$ is again a section. In this way we can now
return to the general, nontrivial bundle situation, formulating
the following


\begin{deff} \label{defbundle}
\noindent A Q-bundle $\pi \colon \cM\to\cM_1$ with typical fiber $\cF$ and
a holonomy algebra $\cG\subset \cD^{<0}(\cF)$ (a chosen graded
Lie subalgebra of vector fields on $\cF$, closed under the derived
bracket)  is a surjective Q-morphism, satisfying the local triviality
 condition: there exists an open cover $\{\cU_i\}$ of $\cM_1$
 such that the restriction of $\pi$ to each $\cU_i$ admits
 a trivialization $\pi^{-1}(\cU_i)\simeq \cU_i\times \cF$ in the
 category of Q-manifolds and this trivialization is
 glued over $\cU_i\cap \cU_j$
 by inner automorphisms which belong to $\exp
(\cG')$ where $\cG'$ is as in (\ref{def_of_G-prime}) with
$\cN =\cU_i\cap \cU_j$.

\noindent A gauge field is a section of $\pi$ in the category of
graded manifolds, that is, a degree preserving map $\gf\colon \cM_1\to\cM$
which obeys $\pi\circ \gf =\mathrm{Id}$. A gauge transformation (an
infinitesimal gauge transformation) is a fiber-wisely acting inner
automorphism (vertical inner derivation) of the total space of
$\pi$.
\end{deff}


\noindent Concatenating a section with a vertical automorphism of the total space,
one obtains an action
of the group of gauge transformations on the space of sections. In generalization of
 (\ref{gauge_vs_derivative}) one then has
 \begin{proposition}\label{variation_of_phi} Given a gauge field $\gf$
 and an infinitesimal gauge transformation $X = ad_Q(Y)$, the
 variation of $\gf$ along $X$ can be identified with the derivation
 $(\de_X \gf)^* := \gf^* X$ covering $\gf$.  \end{proposition}

\section{Field strength and characteristic classes}

\vskip 2mm\npa The obvious role of curvature arises from the fact
that it can be regarded as an obstruction for a map to satisfy the
Maurer-Cartan equation: this motivating example was considered in
the previous section. In the case of general Q-manifolds, the set
of maps between them does not admit a vector space structure any
more, so rather than using the language of Lie brackets, we dealt
with differential graded algebras as a more general
tool. We have also stressed earlier that the Maurer-Cartan (or
zero curvature) equation is a particular example of the chain
property $Q_1 \gf^* -\gf^* Q_2 =0.$
 The operator \cite{BKS}
\begin{equation} F:=Q_1\gf^* -\gf^* Q_2 \, , \label{field_op}
\end{equation}
called {\em the field strength}, being a replacement of the curvature,
is a degree one derivation of functions on the target manifold $\cM_2$
taking values in functions in the source manifold $\cM_1$ and covering
$\gf^*$. (Here $\gf \colon \cM_1 \to \cM_2$ is a morphism of graded
manifolds, corresponding to a gauge field in a trivial bundle. In the
case of a non-trivial bundle, $\cM_2 = \cM$, the total space of a
Q-bundle $\pi \colon \cM \to \cM_1$, and the gauge field $\gf$
satisfies $\pi \circ \gf = \mathrm{Id}$.) We have evidence, for
instance from the Yang-Mills theory, that the curvature is a
meaningful object itself, so one can expect a similar importance of
the ``field strength''.\footnote{This idea was implemented e.g.~in
\cite{Strobl_AYM}.} It is advantageous to reformulate the operator
(\ref{field_op}) somewhat, so that the Leibnitz-property (\ref{Leib})
is replaced by a morphism of algebras (appropriate polynomials should
go into polynomials of field strengths). To this end we regard the
following non-commuting diagram,
\begin{diagram}
T[1]\cM_1 & \rTo^{\gf_*} & T[1]\cM_2\\
\uTo_{Q_1} && \uTo_{Q_2}\\ \cM_1 &\rTo_\gf & \cM_2
\end{diagram}
 where the
homological vector fields are considered as maps and, being of degree
one, the tangent bundle was shifted in degree so that the maps are
morphisms of graded manifolds. Now one notes that both ways from
$\cM_1$ to $T[1]\cM_2$ end in the same fiber over $\cM_2$; thus it is
meaningful to define the difference $f \colon \gf_* \circ Q_1 - Q_2
\circ \gf$, covering $\gf$
\begin{diagram} & & T[1]\cM_2 \\
& \ruTo^f & \dTo\\ \cM_1 &\rTo_\gf & \cM_2
\end{diagram}
It is easy to convince oneself that for any function $h \in
C^\infty(\cM_2)$ and
any $\alpha, \beta \in C^\infty(T[1]\cM_2)$ one has
\begin{equation}
f^* (h) = \gf^* (h) \; , \quad  f^* (\md h) = F
(h) \; , \quad f^*(\alpha \beta) = f^*(\alpha) f^*(\beta) \, .
\label{fstar}
\end{equation}
We shall see below that $f$ is a Q-morphism, if $T[1]\cM_2$ is
endowed with a suitable Q-structure.

\vskip2mm\npa For any graded manifold $\cM$, the algebra of
functions on $T[1]\cM$ admits a simple description as the algebra of
super differential forms $\Omega(\cM)$ (according to the
Bernstein-Leites sign convention \cite{Leites}). More precisely, the algebra of
forms is generated by $h$ and $\md h$ for all functions $h$ with the
following relations: \beq h\;\md h' =(-1)^{\deg (h)(\deg (h')+1) }\md
h' \;h\;, \hspace{3mm} \md (hh')=\md h\; h' +(-1)^{\deg (h) }h\,\md h'
\;.\eeq This algebra is naturally bi-graded by degrees of functions
and orders of forms, such that $\md$, the super-version of the de Rham
differential, becomes a (nilpotent) operator of degree zero with
respect to the first grading and of degree one with respect to the
second grading. The super-commutativity relations are subordinated to
the total grading which is the sum of the two. A vector field $X$ of
degree $p$ gives a contraction of degree $p-1$ acting as follows: \beq
\iota_X \left( f\;\md h\right)=(-1)^{\deg (f)(\deg (X)+1)}f \; X(h)
\;.\eeq The super Lie derivative along $X$, an operator of degree p,
is defined as the commutator \beq\label{Cartan_super} \cL_X :=
\iota_X\md +(-1)^{\deg (X)}\md \iota_X\;.\eeq By construction, $\cL_X$
super-commutes with the de Rham differential and agrees with the
action of vector fields on functions, $\cL_X (f)=X(f)$.  Furthermore,
one can also check that the Lie derivative respects the super-Lie
algebra of vector fields, generalizing the formulas for even
manifolds, such that the following identities hold:
\beq\label{Lie_derivative_identities} [\cL_X,
\cL_Y]=\cL_{[X,Y]}\;, \hspace{3mm} [\cL_X, \iota_Y]=\iota_{[X,Y]}\;.\eeq
In particular, if $Q$ is a homological vector field, we immediately
obtain that\footnote{In Eq.~(\ref{Lie_derivative_identities}) the
brackets indicate \emph{graded} commutators. For an odd vector field
$Q$ the following equations are thus \emph{anti}commutators, the de
Rham differential $\md$ anticommutes with the Lie derivative $\cL_Q$.}
\beq \label{commute} [\md, \cL_Q]=[\cL_Q, \cL_Q]=0\;.\eeq
As a corollary we conclude that the total space of $T[1]\cM$ for a
Q-manifold $\cM$ is a bi-graded manifold supplied with a couple of
super-commuting Q-structures which are of degree one w.r.t.~the
first and the second gradings, respectively. Let us denote the total
differential as $Q_{\scriptscriptstyle T\cM} =\md +\cL_Q$.

\begin{proposition}\label{f-chain_map}
The map $f\colon  \cM_1\to T[1]\cM_2$ is a Q-morphism w.r.t.~the total
Q-structure on the target, that is, the following chain property
holds: \beq Q_1 f^* -f^* Q_{\scriptscriptstyle T\cM_2}=0\;.\eeq
\end{proposition}

\proof Taking into account that the l.h.s.~of the equation is
always a derivation, it is sufficient to apply it on generators of the
algebra of forms (on functions and exact 1-forms). Using
Eqs.~(\ref{fstar}) it is then easy to complete the proof. $\square$

\vskip 2mm\npa There exists also a conceptually more enlightening proof of the
previous proposition. In fact it is easy to convince oneself that
except for $\gf$ all the arrows in the first diagram above are
morphisms of Q-manifolds, if we equip the tangent bundles with the
respective de Rham differentials. Now the map $f \colon \cM_1 \to
T[1]\cM_2$ differs from the (in this sense) Q-morhpism $\gf_* \circ
Q_1$ by a substraction along the tangent fibers by the respective
``value'' of $Q_2$. This substraction corresponds to
$\exp(\iota_{Q_2})$, a diffeomorphism of $\cM_2$ generated by the (in
$T[1]\cM_2$ vertical) vector field $\iota_{Q_2}$, as one can most
easily verify on local coordinate functions. Correspondingly, the new
map, which is our $f$, will remain a chain map, if the de Rham
differential is twisted/conjugated by $\exp(\iota_{Q_2})$. We
summarize this in the following
\begin{lemma} \label{lemtwist} The field strength $f \colon \cM_1 \to
T[1]\cM$ of a gauge field $\gf \colon \cM_1 \to \cM$ can
be defined by the formula $f = \exp(\iota_Q) \circ \gf_* \circ
Q_1$. It is a Q-morphism w.r.t.~
\beq Q_{T\cM}\equiv \md +\cL_Q =\exp(\iota_Q)\,\md\,\exp(-\iota_Q)\;,
\label{twist}
\eeq
\end{lemma}
We remark that $\iota_Q$ does not square to zero since $Q$ is odd;
still the last equality follows easily from the general formula
$\exp(A) \, B \, \exp(-A) = \exp(ad_A) \,B$ (valid for operators
$A$ and $B$ that are not both odd) together with
(\ref{Lie_derivative_identities}) and $[Q,Q]=0$.

\vskip 2mm\npa
A natural example of the chain map property of $f$
is provided by the Weil algebra.  It is well-known that, if one has a
graded morphism from $\Lambda (\g^*)$ of a Lie algebra $\g$ to some
differential graded commutative algebra $\cA$, which is not
necessarily a chain map, we can always extend it as a chain map,
acting from the Weil algebra $W(g)=S^\dot (\g^*)\otimes\Lambda (\g^*)$
to $\cA$. The construction is working as follows: given a graded
morphism $\Lambda (\g^*)\to\cA$, we identify it with some $A$ which
belongs to the dg Lie algebra $\cA\otimes\g$, where the differential
and the bracket are extended by linearity: \beqn \md (\a\otimes
X):=\md\a\otimes X\,,\hspace{3mm} [\a\otimes X,\b\otimes
Y]:=\a\b\otimes[X,Y]\,\eeq  for any $\a,\b\in\cA$ and $X,Y\in\g$.
Defining $F_A:=\md A+\frac{1}{2}[A,A]$ (we recognize the curvature of
a connection in a trivial bundle as a particular example), the
required map $W(g)\to \cA$ is \beq
\Phi\otimes\omega\mto \Phi (\,\underbrace{F_A,\ldots, F_A}_{q\;
\mathrm{times}}\,)\,\omega (\,\underbrace{A,\ldots, A}_{p\;
\mathrm{times}}\,)\;,\hspace{3mm} \Phi\in S^q (\g^*)\;, \omega\in
\Lambda^p (\g^*)\;.\eeq One can easily check that the grading and
differential in the Weil algebra are chosen in such a way that it
becomes isomorphic to $\Omega (\g [1])$ supplied with the
above total differential. Furthermore, the chain map described above
is nothing but our map $f$, if $\cM_2 =\g [1]$ and $\cA=C^\infty
(\cM_1)$.\footnote{In fact, the statement in Proposition
\ref{f-chain_map} can be easily adapted to the situation where
$C^\infty (\cM_i)$ are replaced by
arbitrary differential graded commutative super algebras and
$\gf$ (or better $\gf^*$) by an arbitrary degree preserving
morphism.}

\vskip 3mm\npa Let $\cG$ be a graded Lie subalgebra of vector
fields of negative degree on a Q-manifold $(\cM,Q)$ which is
closed under the derived bracket.

\begin{deff} \label{basic}
A differential form $\omega\in\Omega (\cM)$ is called
a (generalized) $\cG$-basic form, if $\cL_\e
(\omega)=0=\cL_{ad_Q(\e)} (\omega)$ for each $\e\in\cG$. We denote
the space of $\cG$-basic forms as $\Omega (\cM)_\cG$.
\end{deff}

\noindent  $\Omega(\cM)_\cG$ is a  graded
commutative algebra, which is stable with respect to both
differentials: The product of two such forms is again a
$\cG$-basic form. The Lie derivative commutes with the de Rham
operator, so the space is invariant with respect to $\md$.
Furthermore, by the identities (\ref{Lie_derivative_identities}),
one has $[\cL_Q, \cL_\e]=\cL_{ad_Q(\e)}$ and
$[\cL_Q,\cL_{ad_Q(\e)} ]=0$, thus $\Omega(\cM)_\cG$ is also
closed with respect to $\cL_Q$.

\vskip 3mm\npa Now we apply this machinery in a rather
straightforward way for the construction of characteristic classes
associated to any section (gauge field) of a Q-bundle.

\begin{theorem}\label{construction_of_charclasses} Let $\pi
\colon \cM\to\cN$ be a Q-bundle with a typical fiber $\cF$, a  holonomy
algebra $\cG$, and $\gf$ a section of $\pi$ (in the graded
sense)---cf.~Definition \ref{defbundle}. Then there is a well-defined
map in cohomology
\beq\label{characteristic_map} H^p (\Omega (\cF)_\cG,
Q_{\scriptscriptstyle T\cF})\to H^p (C^\infty (\cN),
Q_{\scriptscriptstyle \cN} )\;,\eeq which does not depend on
homotopies of $\gf$.
\end{theorem}

\begin{lemma}\label{variation_of_f}
 Let $(\cM_1,Q_1)$ and $(\cM_2,Q_2)$ be Q-manifolds, $\gf$ a morphism
 $\cM_1\to\cM_2$ of the underlying graded manifolds, and $Y$ a vector
 field on $\cM_2$ of degree minus one, generating the inner derivation
$X =ad_{Q_2} (Y)$. Then the induced variation---the induced
infinitesimal gauge transformation---of $f^* \colon \Omega(\cM_2) \to
 C^\infty(\cM_1)$ is given by:
 $\de_{X} f^* =f^*\cL_X $.
\end{lemma}
\proof Since $f^*$ is a degree preserving map, its variation
w.r.t.~an infinitesimal flow, generated by a degree zero vector field,
is a degree preserving derivation above $f^*$. Thus it is sufficient
to check the identity on functions $h$ and exact 1-forms $\md h$ over
$\cM_2$.  For the first part we can use Proposition
\ref{variation_of_phi}, since $f^* (h)=a^* (h)$. Likewisely, using this
proposition and eqs.~(\ref{fstar}),
 we find: \beqn \de_X f^* \md
h=\de_X F (h)=Q_1\de_X\gf^* (h)- \de_X\gf^*Q_2 (h)=Q_1\gf^*
X (h)- \gf^* X Q_2 (h)\;.\eeq But, as an inner derivation, $X=ad_{Q_2} (Y)$
commutes with $Q_2$. Thus, $ \de_X f^* (\md
h)= F\cL_X(f) = f^* \cL_X( \md f)$, which concludes the proof.
$\square$

\begin{lemma}\label{gauge_invariance}
Let $\omega\in\Omega (\cM_2)_\cG$, $\cG$ vector fields on $\cM_2$
closed w.r.t.~the Lie and the derived bracket,  and $\omega'$ be its trivial
extension to the total space of $\pi \colon \cM_1\times\cM_2\to\cM_1$.
Then $f^* (\omega ')$ is invariant under the action of $\exp (\cG')$ on sections,
where $ \cG' :=ad_{Q}
\left(C^\infty (\cM_1, \cG)\right)\cap \cD^0 (\cM_1 \times \cM_2)$.
($\gf \colon \cM_1 \to \cM_1 \times \cM_2$,  $Q\equiv Q_1+Q_2$).
\end{lemma}

\proof Let us replace the target
manifold $\cM_2$ in Lemma \ref{variation_of_f} by the total space $\cM
= \cM_1 \times \cM_2$ and an arbitrary graded morphism $\gf$ with a
section of $\pi$. Suppose we are given an infinitesimal gauge
transformation $X =ad_Q (Y)$. Applying Lemma
\ref{variation_of_f} to the variation of $f^*$ along $X$ and formula
(\ref{ad_on_total_space}), we obtain:
\beq\label{variation_f_omega} \de_X f^*\omega'=f^*\cL_{ad_Q
(Y)}\omega'=\sum_j f^*\left(\cL_{Q_1 (b^j)Y_j}+ (-1)^{\deg
(b^j)}\cL_{b^j ad_{Q_2} (Y_j)}\right)\omega'\;.\eeq We are left to
prove that $f^*\cL_{hv}\omega '=0$ for each $h\in C^\infty
(\cM_1)$ of degree q and $v\in\cD^p (\cM_2)$ which obeys $\cL_v
\omega' =0$: By the definition of Lie derivative, \beqn
\cL_{hv}\,\omega'=\left( \iota_{hv}\md +(-1)^{p+q}\md
\iota_{hv}\right)\omega' = \left((-1)^{p+q}\md h\wedge
\iota_v+h\cL_{v}\right)\omega'=(-1)^{p+q}\md h\wedge \iota_v\omega'\;.\eeq
On the other hand, for any section $\gf$ its field strength $F$
is a vertical derivation:  \beqn F (h)=Q_1\gf^* \pi^*
(h)-\gf^*(Q_1+Q_2)\pi^* (h)=0\;,\hspace{3mm} \forall h\in
C^\infty (\cM_1)\,,\eeq since $\pi \gf\equiv \mathrm{Id}$ and
$Q\pi^*=\pi^* Q_1$. With this equation and (\ref{fstar}) we then indeed obtain \beqn
f^*\cL_{hv}\,\omega'=(-1)^{p+q} f^* \left( \md h\wedge
\iota_v\omega'\right)=(-1)^{p+q} F(h)f^* (\iota_v\omega')\equiv 0\;.\eeq$\square$

\vskip 4mm\noindent{\bf Proof of Theorem
\ref{construction_of_charclasses}~}

\vskip 2mm\noindent Let us fix a trivialization of $\pi$ over an
open cover $\cU_i$. Suppose we are given a section $\gf$, then
there is a family of sections $\gf_i$ over $\cU_i$ which are the
restrictions of $\gf$. Applying $f_i^*$ to each trivial extension
$\omega'_i$ of $\omega$ to $\cU_i\times\cF$, we obtain a family of
$Q_1$-cocycles in $C^\infty (\cU_i)$, denoted as $\mathrm{char}_i
(\omega)$. Taking into account that the trivialization is glued by
a transition cocycle of gauge transformations belonging to
$\exp (\cG')$---cf.~definition \ref{defbundle}---,
$\gf_i$ and $\gf_j$ are related by a gauge
transformation; thus, applying Lemma \ref{variation_f_omega}, we
obtain $\mathrm{char}_i (\omega)=\mathrm{char}_j (\omega)$
over $\cU_i\cap\cU_j$. Thus one has a global $Q_1-$cocycle
$\mathrm{char} (\omega)$, such that $\mathrm{char}_i (\omega)$ is
its restriction to $\cU_i$.

\vskip 2mm\noindent It remains to prove that if $\gf(t)$, $t\in
[0,1]$, is a smooth family of sections, then the cohomology class of
the corresponding $\mathrm{char} (\omega)(t)$ does not change. For
this purpose we use the same argument as for the usual Chern-Weil
formalism. A family of sections $\gf(t)$ can be thought of as a
section $\tilde \gf$ of the following extension of $\pi$: \beqn
\tilde\pi \colon \cM\times T[1]I\to \cN\times T[1]I\,,\hspace{3mm}
I=[0,1]\,.\eeq By construction, for any function $h$ on the total
space of $\tilde\pi$, written as $h=h_0 +\md t \,h_1$, where $h_i \in
C^\infty (\cM\times I)$, the pull-back with respect to $\tilde \gf$ is
$\tilde \gf^* h =\gf^* (t)h_0+\md t\,\gf^* (t)h_1$. The new field
strength operator is $\tilde F =F(t) +\md t\, \pt_t \gf^*
(t)$. Applying the corresponding characteristic map to $\omega$, we
obtain a cocycle on $\cN\times T[1]I$, which decomposes as follows:
\beq
\widetilde{\mathrm{char}} (\omega)=\mathrm{char} (\omega)(t) +\md t\,
\beta \;.\eeq Using the closedness with respect to
$Q_{\scriptscriptstyle \cN}+\md_I$, we immediately get the identity
\beq \pt_t \mathrm{char} (\omega)(t)= Q_{\scriptscriptstyle \cN}\beta
(t)\;,\eeq which implies the invariance of characteristic classes in
cohomology:
\beqn \mathrm{char} (\omega)(1)-\mathrm{char}
(\omega)(0)=Q_{\scriptscriptstyle \cN}\int\limits_0^1 \md t\,
\beta (t)\;.\eeq This completes the proof of Theorem
\ref{construction_of_charclasses}. $\square$

For some purposes like the construction of secondary characteristic classes,
it may be useful to display the
transgression $\beta(t)$ explicitely. If $q^\alpha$ denote (local)
graded coordinates on the fiber $\cF$, and
on the chart $\cU_i$ we use the notation
$\gf_i(t)^*(q^\alpha) =: A^\alpha(t)$ and
$f_i(t)^*(q^\alpha) =: F^\alpha(t)$, denoting the corresponding
tower of gauge fields and field strengths, then for
$\omega = \frac{1}{p!}\omega_{\a_1 \ldots \a_p}(x) \md x^{\a_1} \wedge \ldots \wedge
\md x^{\a_p}$ one finds
\beq \label{beta}
\beta(t)|_{\cU_i} = \frac{1}{(p-1)!} \,
\gf_i(t)^* \left(\omega_{\a_1 \ldots \a_p} \right) \partial_t
\left(A^{\alpha_1}(t)\right) \wedge F^{\a_2}(t) \wedge \ldots \wedge F^{\a_p}(t)
\eeq
Note that if $\cF$ is a Qk-manifold and $\cN=T[1]\Sigma$, then $A^\a$s
are in general a tower of differential forms of degree zero up to degree
$k$ (and likewise $F^\a$s differential forms of degree one up to
degree $k+1$). Also, in general $\omega_{\a_1 \ldots \a_p}$ will depend on
$x^\a$ and thus its pullback by  $\gf_i(t)^*$ produce a function in the
gauge fields. Clearly, by construction, $\beta(t)$ is well-defined globally
on the base $\cN$.

\vskip 2mm\npa By use of Theorem \ref{construction_of_charclasses}
we associate a characteristic class to each gauge field and
cohomology class of the subcomplex $(\Omega (\cF\!)_\cG,
Q_{\scriptscriptstyle T\cF})$ of $\cG-$invariant forms on the
fiber.
\begin{proposition} \label{proptriv} Suppose $c$ belongs to the kernel in cohomology
of the
canonical map $(\Omega (\cF\!)_\cG, Q_{\scriptscriptstyle
T\cF})\to (\Omega (\cF\!), Q_{\scriptscriptstyle T\cF})$. Then the
corresponding characteristic class is trivial for any trivial
Q-bundle $\cM\to\cN$.
\end{proposition}
\proof If the Q-bundle is trivial, then there exists a projection
of the total space to the fiber $p_{\scriptscriptstyle \cF}:
\cM\to\cF$ which is a Q-morphism. The corresponding characteristic
class attached to a gauge field $\gf$, which is nothing but $f^*
p_{\scriptscriptstyle \cF}^*c$, is certainly trivial on the base,
$f$ being a Q-morphism.
  $\square$



\vskip3mm\npa Bringing the gauge transformation of the field strength
$\de_{ad_{Q_2}(Y)} f^*$ in Lemma \ref{variation_of_f} into the form
$f^* \circ ad_{Q_{T\cM_2}}(\cL_Y)$, we observe that it fits the
pattern of Proposition \ref{variation_of_phi}: Given a gauge field
$\gf \colon \cM_1 \to \cM$ in a Q-bundle $\cM \to \cM_1$ and an
infinitesimal gauge transformation generated by $ad_Q(Y)$, $Y$ a
vertical vector field of degree minus one, we can canonically
associate to these data: A Q-bundle $(\widetilde{\cM} \to \cM_1,
\widetilde{Q})$, a gauge field $\widetilde \gf$ which is a true
section of the Q-bundle (i.e.~indeed in the category of Q-manifolds),
and a vertical vector field $\widetilde{Y}$; here $\widetilde{\cM} =
T[1]\cM$, $\widetilde{Q}= Q_{T\cM}
\equiv \md + \cL_Q$, $\widetilde \gf=f$ (a true section due  to Proposition
\ref{f-chain_map}), and  $\widetilde{Y}=\cL_Y$.
Using these identifications we can
always, vice versa, recover the primary data. The extended Q-bundle is
quite special, however: vector fields in the gauge transformations
have a very particular form, generated by
$\widetilde{Y}=[\iota_Y,\md]$ for some $Y$ living on $\cM$, and
likewisely restricted is the holonomy group $\widetilde{\cG} \cong
\cG$ of $\widetilde{\cM}$. Still, in this language, a $\cG$-basic form on $\cF$
translates into a function on the new fiber $\widetilde{\cF} =
T[1]\cF$ which is annihilated by $\widetilde{\epsilon}$ and
$[Q_{\widetilde{\cF}},\widetilde{\epsilon}]$ where again
$\widetilde{\epsilon}$ is of the particular form
$[\iota_\epsilon,\md]$ (for some $\epsilon$ living on $\cF$ generating
the holonomy $\cG$ as before). 

\section{Some applications}

\npa The first evident example is a principal G-bundle over $M$.
The corresponding Q-bundle is provided by the anchor map $\rho\colon
E[1]\to T[1]M$ of the associated Atiyah algebroid $E$, cf.~examples
(3) and (4) in paragraph 1.3. The holonomy algebra is $\g$, the Lie
algebra of $G$, and the fiber is $\g [1]$. Using the isomorphism
$W(\g)\simeq \Omega (\g [1])$ explained in the previous section, we
can easily see the isomorphism between the space of basic forms and
symmetric $G-$invariant polynomials on $\g^*$, $\Omega (\g [1])_\g
\simeq S(\g^*)^G$: An element of $\omega \in \Omega (\g
[1])_\g$ always has the form $\omega = \frac{1}{p!}
\omega_{a_1 \ldots a_p} \md
\xi^{a_1} \wedge \ldots \wedge
\md \xi^{a_p}$, if $\xi^a$ denote the odd coordinates on $\g[1]$, with
$\omega_{a_1 \ldots a_p}$ being constant, completely symmetric, and
adinvariant (following from invariance w.r.t.~$\cL_\e$, graded
antisymmetry of $\omega$, and invariance w.r.t.~$\cL_{ad_Q(\e)}$,
respectively).

\vskip2mm \noindent As explained before, a section of $\rho$ is  a
connection in the principal $G-$bundle, and the construction of
characteristic classes in Theorem \ref{construction_of_charclasses}
applied to this particular case reproduces the Chern-Weil map. Indeed,
with $F^a$ denoting the local curvature 2-forms, from the above
$\omega$ one obtains $$ \frac{1}{p!}  \omega_{a_1 \ldots a_p} F^{a_1}
\wedge \ldots \wedge F^{a_p} \, . $$ Likewise, $\b(t)$ of
Eq.~(\ref{beta}) specified to this case, gives the standard
transgression formula in this example.

\vskip 3mm\npa Equivariant cohomology and gauging of WZ-terms in
sigma models: A near-at-hand extension of the Weil algebra $W(\g)$
is the Weil model of equivariant cohomology. Let $G$ be a Lie
group acting on a manifold $M$ and $\rho \colon \g \to \cD(M)$,
$\cD(M) \equiv \Gamma(TM)$, be the corresponding Lie algebra
action. The complex one looks at is $C=W(\g) \otimes \Omega(M)$,
equipped with the sum of the previously introduced differential on
$W(\g)$ and the de Rham differential on the forms on $M$. In
Q-language this is the space of functions on
$\widetilde{\cM_2}:=T[1] (\g[1]\times M)$ and the differential
gives rise to a homological vector field on it, which we want to
call $Q_W$. The $\g$-action is extended in a natural way to this
complex. To describe this in the language of section
\ref{sectionQbundles}, we need a map from $\g$ into vector fields
$\cD(\widetilde{\cM_2})\ni X$ of degree minus one such that
$ad_{Q_W}(X)$ generates the $\g$ action. This is easy to find:
Take an element $\e$ in the Lie algebra, it generates canonically
a vector field of degree minus one on $\g[1]$ and thus (by lifting
as a Lie derivative) also on $T[1]\g[1]$. Likewise $\rho(\e)$
gives a vector field on $M$, its contraction $\iota_{\rho(\e)}$ is
a degree minus one vector field on $T[1]M$. We denote the sum of
these two vector fields by $i_\e$. Now $ad_{Q_W}(i_\e) \equiv
[Q_W,i_\e]$ is easily verified to generate the canonical diagonal
$\g$-action on $C^\infty(\cM_2)$.

An element $\alpha \in C$ is called horizontal if, in the above
language, it is annihilated by $i_\e$ for all $\e \in \g$. If, in
addition, it is also $\g$-invariant, i.e.~also annihilated by
$ad_{Q_W}(i_\e)$, it is called basic. The space of basic elements in
$C$ is denoted by $C_\g$ or by $\Omega(M)_G$. Although the space
$\widetilde{\cM_2}$ is of the form $T[1](\cM_2)$ with $\cM_2 \equiv
\g[1] \times M$ \emph{and} the homological vector field $Q_W$ is nothing but
the total differential of the Chevalley-Eilenberg differential $
Q_{CE} \simeq \md_\g$ on $\g[1]$ extended trivially to $\cM_2$, $Q_W =
\md + \cL_{Q_{CE}}$, the notion of basic elements does \emph{not}
agree with the one of Defintion \ref{basic}. This, as we will see also
in more detail below, is related to the last remark in the previous
section; $\widetilde \e \equiv i_\e= \cL_\e + \iota_{\rho(\e)}$ is of
the required form of a Lie derivative only on $T[1]\g[1]$, but not
also on $T[1]M$.

The reason for this apparent discrepancy is that not all the data
have been incorporated properly into the $Q$-structure on $\cM_2$;
in particular, the representation $\rho$ entered only when
considering $\widetilde \e$. This can be cured easily, however,
and will lead us automatically to the so-called Cartan model of
equivariant cohomology. The initial data give rise to an action
Lie algebroid $E = M \times \g$ over $M$, $\rho$ yielding its
anchor map. The respective homological vector field of $\cM_2 :=
E[1]$ has the form \beq \label{actionLie} Q=\rho + Q_{CE} \eeq
where we interpreted $\rho$ as an element in $\g^* \times \cD(M)$,
viewing $\g^*$ as a linear and thus degree one function on
$\g[1]$. For non-abelian $\g$, the vector field $\rho$ does not
square to zero, but it is easily verified that $\rho^2 = -
[Q_{CE},\rho]$, so that indeed $Q^2 = 0$.\footnote{In local
coordinates $x^i,\xi^a$ on $\cM_2=M\times \g[1]$, one recognizes
in $Q= \xi^a \rho_a  - \frac{1}{2} \xi^b \xi^c C^a_{bc}
\partial_a$, where $C^a_{bc}$ are the structure constants of $\g$
in the basis $\xi_a$ dual to $\xi^a$ and $\rho_a \equiv
\rho(\xi_a)$, one recognizes the standard YM-type BRST charge.
Although not inspired by  \cite{Kalkman}, the considerations in
this paragraph partially parallel, and possibly also simplify and
highlight, those of that paper.} Let us call the canonical lift
$\md + \cL_Q$ of $Q$ to $T[1]\cM_2$ by $Q_C$. Since $Q$ differs
from $Q_{CE}$ by the addition of $\rho$ and vertical vector fields
on $T[1]\cM_2\to \cM_2$ being contractions of vector fields coming
from the base always (super)commute, it follows immediately from
(\ref{actionLie}) and Lemma \ref{lemtwist} that \beq Q_W =
\exp(-\iota_\rho)\: Q_C \: \exp(\iota_\rho) \, . \eeq In this more
geometric picture we find that indeed the notion of basic above
agrees with Defintion \ref{basic}: $\epsilon$ as before, we see
that
\begin{eqnarray}  \exp(-\iota_\rho) \: \cL_\e \: \exp(\iota_\rho) &=&
i_\e \\  \exp(-\iota_\rho) \:
\cL_{ad_{Q}(\e)} \: \exp(\iota_\rho) &=& ad_{Q_W}(i_\e) \label{gen}
\end{eqnarray}
where for the first equality we made use of $[\cL_\e,\i_\rho] =
\iota_{[\e,\rho]} = \iota_{\rho(\e)}$, which in turn commutes with
$\iota_\rho$, from which the second one follows immediately on
observing $\cL_{ad_{Q}(\e)}=[Q_C,\cL_\e]$. Moreover, elements in
$C \cong \Omega(E[1])$ annihilated by $\cL_e$ are easily seen to
be elements in $S^\dot (\g^*) \otimes \Omega(M)$. On the other
hand, the map $\rho$ is equivariant and thus the respective
homological vector field $\rho$ is $\g$-invariant. This implies
that $ad_{Q_W}(i_\e)$ commutes with $\exp(\iota_\rho)$, which in
addition to eq.~(\ref{gen}) yields the equality $\cL_{ad_{Q}(\e)}
= ad_{Q_W}(i_\e)$. Thus annihilation by $\cL_{ad_{Q}(\e)}$ implies
$\g$-invariance and basic differential forms on the action Lie
algebroid $E[1]$ following Definition \ref{basic} correspond
precisely to elements in $\left(S^\dot (\g^*) \otimes
\Omega(M)\right)^G$, which is the Cartan model of equivariant
cohomology.

This is now also the right language and departure point for the
discussion of gauging of WZ-terms in sigma models (cf.~also
\cite{Figueroa-Stanciu}). Let d-dimensional ``spacetime'' be the
boundary of $N$ and consider the space of maps $X \colon N \to M$
as (part of) the ``fields'' of the sigma model. A WZ-term then is
induced by a closed (d+1)-form $H$ on $M$, $S_{WZ}[X] = \int_N
X^*H$. G-invariance of $H$ yields $S_{WZ}$ invariant under
``rigid'' G-transformations, but not under ``local'' ones,
i.e.~where the transformation parameters are permitted to vary
along $N$. To capture this fact in the present framework, we
extend $X$ trivially to a map $\gf_0 \colon (T[1]N,Q_{dR}) \to
(M,0)$; the WZ-term then can also be written as
$S_{WZ}[\gf_0]=\int_N f_0^* H$, where $f_0=X_* \colon T[1]N\to
T[1]M$ is the ``field strength'' of $\gf_0$. In this simple case
the chain map property of Proposition \ref{f-chain_map} reduces to
the wellknown fact that the pullback map $X^*$ commutes with the
de Rham differential.

The above mentioned ``local G-transformations'' will now become gauge
transformations in the trivial Q-bundle $T[1]N \times M \to T[1]N$ as
discussed in section \ref{sectionQbundles}. Indeed, the representation
$\rho$ singles out a Lie subalgebra $\g \subset \cD(M)$ of symmetries
on $M$ (there are certainy no inner infinitesimal symmetries on
$(M,0)$); thus functions on $T[1]N$ with values in this Lie subalgebra
of degree zero vector fields on $M$ can be used as a proper
replacement of the infinitesimal gauge Lie algebra $\cong \cG'\ni X$
in this context. Its lift to the field strength, $\delta_X f_0^* =
f_0^* \cL_X$ (cf.~Lemma \ref{variation_of_f}), produces the correct
transformation. Note that now $\gf_0$ is considered as a section in
the above trivial Q-bundle (and $H$ is extended in a likewise manner
to the trivial bundle); only like this \cite{BKS} we can accomodate
for the $N$-dependence of the infinitesimal generator $X$ of the
transformations, $X=\e^a \rho_a$, in the notation of the previous
footnote, with $\e^a$ an arbitrary function on $N$. Clearly,
$S_{WZ}[\gf_0]$ is \emph{not} invariant under any such gauge
transformation, since the Lie derivative contains the de Rham
differential on $N$ and one obtains $f_0^* \cL_X H = \md \e^a \wedge
X^* (\iota_{\rho_a} H) + \e^a X^*( \cL_{\rho_a} H)$, where only the
second term vanishes by invariance of $H$.\footnote{If one considers
$f_0$ as a ``gauge field'' $\widetilde{\gf_0}$ itself, cf.~the
discussion at the end of the previous section, the gauge
transformations become inner and this becomes a special case of
eq.~(\ref{ad_on_total_space}). However, the field $\widetilde{\gf_0}$
is restricte to derive from the ``field strength'' of some $\gf_0$ and
it is also this perspective that now permits to discuss the gauging of
the WZ term in a concise manner.}

To cure this one wants to introduce extra gauge field dependent terms
into the action functional, i.e.~terms depending on additional
$\g$-valued 1-forms $A^a$ on (the boundary of) $N$. This is possible
if $H$ permits a G-equivariantly closed extension $\hat H$, cf.~also
\cite{Figueroa-Stanciu}. In our picture the resulting
invariant action functional is now easy to obtain: We simply replace
$(M,0)$ in all of the constructions above by the action Lie algebroid
$E[1]=M \times \g[1]$ with its canonical differential,
eq.~(\ref{actionLie}). $\gf$ now is a section in the trivial Q-bundle
$T[1]N \times E[1] \to T[1]N$ and $\hat H$ extended trivially from
$\Omega(E[1])$ to a differential form on the total Q-bundle (analogous
to the case of $H$ before). Gauge transformations are now inner right
away, they are generated by $\e$s as described in the above
Cartan-type model of equivariant cohomology, tensored with functions
on the base. Gauge invariance of $f^* \hat H$ now follows at once from
the general result of Lemma \ref{gauge_invariance}. From Proposition
\ref{f-chain_map}, moreover, it follows immediately that
$f^* \hat H$ is closed and that $f^* \hat H-f^*H$ is exact; thus, the
additional gauge fields, corresponding to a degree preserving map from
$T[1]N$ to $\g[1]$, indeed need to be defined over the boundary of $N$
only.


The formalism developed in this paper is certainly aimed at also more
general type of gauge theories as those stemming from a structural Lie
group like in this paragraph. We intend to make this kind of application
more explicit elsewhere, focusing in the present draft mainly on the
issue of (generalized) characteristic classes.\footnote{But cf.~also
\cite{BKS,Strobl_AYM,Gruetzmann-Strobl}, as well as the following paragraph, used to prepare
grounds for characteristic classes associated to ``PQ''-bundles in
the subsequent paragraph 4.7 below.}

\vskip 3mm\npa Another application of the considerations of the previous
section is the following
\begin{theorem} Let $(\cS,\omega)$ be a symplectic Qp-manifold, $p\in \mathbb{N+}$,
as in example (3) of paragraph 1.2, $N$ a $(p+2)$-dimensional
manifold with boundary $\partial N=\Sigma$, and $\gf$ a (degree
preserving) map from  $T[1]N$ to $\cS$. Then
\beq 
\int_{N} f^* \omega =  S^{AKSZ}_{\S,(cl)}\eeq
where $S^{AKSZ}_{\S,(cl)}$  is the (classical part of the)
topological sigma model on the $(p+1)$-dimensional
$\Sigma$ obtained by the AKSZ-method \cite{AKSZ}. \label{theoAKSZ}
\end{theorem}
If $\gf$ is a gauge field in the sense of this paper, i.e., being a
degree preserving map, $\gf$ has degree zero, one obtains the classical
part of the topological action. Permitting all possible degrees of
$\gf$, one gets its BV extension, satisfying the classical BV-master
equation, i.e.~squaring to zero w.r.t.~the BV bracket. For $p=1$ the
action reproduces the Poisson sigma model
\cite{Schaller-Strobl,Ikeda}, for $p=2$ one obtains the Courant sigma
model \cite{Ikeda_3d,Park,Roytenberg06}; the formula above holds for
arbitrary dimensions.

\vskip3mm \noindent Before proving this, we make some general remarks
on PQ-manifolds with $p>0$\footnote{This assumption will be kept without
further mention.}. This is a Q-manifold $\cS$ equipped with a
compatible symplectic form $\omega$ of degree p, i.e.~it obeys
$\cL_\xi
\omega =p\omega$, where $\xi$ is the Euler field which provides the
$\Z -$grading on $\cS$. Non-degeneracy of $\omega$ implies that $\cS$
has degree at most $p$ (if $\cS$ has a nontrivial body, i.e.~its
algebra of functions has degree zero elements, this bound is also
necessarily saturated); since a lower-degree Q-manifold can also be
considered as a degenerate Qp-manifold for some higher p, we will
follow the convention that PQ-manifolds of degree p imply that
$\omega$ has degree p. E.g.~given a quadratic Lie algebra
$(\g,\kappa)$, $\kappa$ denoting the adinvariant scalar product, we
will view $(\g[1],\omega)$, $\omega \sim \kappa$, as a degree 2
PQ-manifold. For $p > 0$ the symplectic form is necessarily exact: By
(\ref{Cartan_super}), one has $\md \iota_\xi\omega =p\omega$, thus
$\omega =\md\a$, where $\a =\frac{1}{p}\iota_\xi\omega$. For any
function $h$ of degree $q$ on $\cS$ we define its Hamiltonian vector
field $X_h$ of degree $q-p$ by the formula: $\iota_{X_h}\omega
=(-1)^{q+1}\md h$. Then Hamiltonian vector fields satisfy the known
relations from the ungraded case:
\beq\label{super_hamiltonian_commutator}[X_{h_1},
X_{h_2}]=X_{\{h_1, h_2\}}\,,\eeq where $\{\cdot , \cdot \}$ is the
induced Poisson bracket of degree $-p$. Now, using the relations
(\ref{Lie_derivative_identities}), we can easily verify that
compatibility of $\omega$ with $Q$, $\cL_Q \omega =0$, implies that
$Q$ is always Hamiltonian: $i_Q\omega =(-1)^p \md\cQ$ with the
Hamiltonian function $\cQ$ of degree p+1.  We summarize this in the
following
\begin{lemma} \label{exact}
For a Qp-manifold with compatible sympectic form $\omega$ one has
\beq  \omega =\md\a \, , \; \a \equiv \frac{1}{p}\iota_\xi\omega \, , \qquad
Q = X_\cQ \, , \; \cQ \equiv \frac{p}{p+1} (-1)^p \iota_Q \alpha \, , \label{Ham}
\eeq
where $\xi$ is the Euler vector field on $\cS$ and the Hamiltonian $\cQ$ of the homological
vector field $Q$ satisfies the master equation  $\{\cQ,
\cQ\}=0$.
\end{lemma}
The last statement follows from eq.~(\ref{super_hamiltonian_commutator}) and the fact that
in positive degrees the only constant is zero ($\cQ$ has degree $p+1$ and thus $\{\cQ,
\cQ\}$ degree $p+2$). From the above
we derive
\begin{lemma} \label{propos_transgr_for_AKSZ} The following
transgression formula holds: \beq\label{formula_transgr_AKSZ}
\omega =Q_{\scriptscriptstyle T\cS} \left(\hat \a\right) \, , \qquad \hat \a =
\a + \frac{(-1)^p}{p}\cQ \equiv
 \frac{1}{p}\left(1 + \frac{1}{p+1} \iota_Q\right)\iota_\xi\omega
\;.\eeq
\end{lemma}

\proof With Lemma \ref{exact} and eq.~(\ref{twist}) we have $\exp(\iota_Q) \omega =
Q_{\scriptscriptstyle T\cS}  \exp(\iota_Q) \a$. The l.h.s.~of this equation
is $\omega + (-1)^p \md \cQ$,
 by defintion of the Hamiltonian for $Q$ and the fact that
$\iota_Q \iota_Q \omega$ vanishes on behalf of  $\{\cQ, \cQ\}=0$. The master equation also
implies that $\md \cQ=Q_{\scriptscriptstyle T\cS} \cQ$, from which we now can derive easily
the wished for transgression formula.  $\square$

\vskip 4mm\noindent{\bf Proof of Theorem
\ref{theoAKSZ}~}

\vskip 2mm\noindent With Lemma \ref{propos_transgr_for_AKSZ} one
obtains $f^* \omega = \md f^* \hat \a$, where the chain property
Prop.~\ref{f-chain_map} of $f$ has been used and the fact that $Q_1$
is just the de Rham differential on $N$ here. Using Stokes theorem,
we are thus left with showing that $f^* \hat \a$ indeed agrees with
the AKSZ action (where we can replace $N$ by its boundary $\Sigma$ now).
By means of the formulas (\ref{fstar}), one can
convince oneself that for any 1-form $\a$ on $\cS$ one has \beq
f^*\a =  \left(\iota_{Q_{T[1]\S}} \gf^* - \gf^* \iota_{Q_\cS}\right) \a \, .
\eeq
Together with the first equation in  (\ref{fstar}) and the second equation in (\ref{Ham}),
this implies
\beq f^* \hat \a =
\iota_{\md_\S}\gf^*\a + (-1)^{p+1}\gf^* (\cQ ) \;, \eeq
where  $Q_{T[1]\S}=\md_\S$ has been used. This expression
agrees precisely with the one found in
\cite{Roytenberg06} for the AKSZ sigma model, which thus
completes the proof. $\square$

\vskip 3mm\npa Q-bundles with PQ-manifolds $(\cS,\omega)$ as
fibers are natural candidates for
a non-trivial characterisitic class along the lines of Theorem
\ref{construction_of_charclasses}. As we saw above, the symplectic form
 $\omega$ itself is closed w.r.t.~the total differential
 $Q_{\scriptscriptstyle T\cS}$ on the fiber and for the holonomy group
 $\cG$ the Lie algebra of (all or a closed subset of the) Hamiltonian
 vector fields of negative degree lends itself naturally, since
 $\omega$ is then also basic w.r.t.~$\cG$.\footnote{In fact, $\omega$
 is even exact within the complex $(\Omega(\cS),Q_{T\cS})$, cf.~Lemma
 \ref{propos_transgr_for_AKSZ} above; however, in general it will fail
 to be exact within the restricted complex
 $(\Omega(\cS)_\cG,Q_{T\cS})$ of basic forms. This happens already for
 the standard characteristic classes, paragraph 4.1 above, where the
 unrestricted cohomology, being isomorphic to deRham cohomology on the
 Lie algebra, cf.~Lemma \ref{lemtwist}, is obviously trivial. It is
 also in this context where Proposition \ref{proptriv} comes into
 play.}  A PQ-bundle (a Q-bundle with PQ-fibers) thus carries a canonical
characteristic class. For the Atiyah algebroid of a principal $G-$bundle
(cf.~example (3) and (4) in paragraph 1.3)
where the Lie algebra $\g$ of $G$ is equipped with a non-degenerate
invariant symmetric form, the corresponding PQ-bundle has a
typical fiber $\g[1]$ together with a degree $p=2$ symplectic form $\omega$
provided by the invariant metric. The canonical characteristic
class is nothing but the second Chern class (or the first Pontrjagin class) of the principal
$G-$bundle and the Theorem \ref{theoAKSZ} simply
gives the well-known local statement: $``\mathrm{Second\;
Chern\; form}=\md \left(\mathrm{Chern\!-\!Simons\;
form}\right)''$. For an arbitrary PQ-bundle over a base $T[1]N$
for a smooth manifold $N$ one thus has a
straightforward generalization of the second Chern class, which is a $p+2$
cohomology class for a degree $p$ QP-fiber.

\vskip 2mm\noindent For example, in the case of $p=1$ a typical fiber $\cF$ is
necessarily of the form $T^*[1]P$ for some Poisson manifold $(P,\{
\cdot , \cdot \})$. Like any Q1-manifold, the total space corresponds
to a Lie algebroid $E$ living over some base manifold $M$. Thus the
Q-bundle has the form
\begin{diagram}
E[1] & \rTo_{\pi} & T[1]N\\
\dTo && \dTo\\ M &\rTo_{\pi_0} & N
\end{diagram}
covering an ordinary bundle $\pi_0 \colon M \to N$ whose typical
fiber is $P$. It can be shown that this bundle $\pi_0$ is a
Hamiltonian Poisson fibration, i.e.~the Lie algebra of its
holonomy group consist of Hamiltonian vector fields of the Poisson
manifold $P$. As we will show in detail in a separate note
\cite{Kotov-Strobl-inprep}, the 3-form class on $N$ that one
obtains in this way does not depend on the chosen gauge field
$\varphi\colon T[1]N \to E[1]$. Moreover, it constitutes an
obstruction to lifting the bundle $\pi_0$ to one where the
transition cocycle takes values in a group whose Lie algebra is
$(C^\infty(P),\{ \cdot , \cdot \})$---while the existence of a
lift to the case of Hamiltonian functions \emph{modulo} constants
is already guaranteed by the existence of the Q-bundle $\pi$
covering $\pi_0$.

\vskip 3mm\npa Suppose we are given a bundle $p\colon \cM\to\cN$
in the category of Q-manifolds which, in general, is \emph{not} locally
trivial. Let us denote by $\Omega_+ (\cN)$ the ideal of all
differential forms the order of which as a differential
form is greater or
equal to one, and by $\cI$ the ideal in $\Omega (\cM)$ generated by
the pullback $p^* \Omega_+ (\cN)$. It is clear that $\cI$ is
closed with respect to both differentials $\md$ and $\cL_Q$, and
thus with respect to the total differential $Q_{\scriptscriptstyle
T\cM}$. For each gauge field $\gf\colon \cN\to\cM$ the
corresponding map $f^*\colon \Omega (\cM)\to C^\infty (\cN)$
vanishes on $\cI$ (this fact was used in the proof of Lemma
\ref{gauge_invariance}). Therefore we have a well-defined chain
map of complexes \beq\label{I-chain} (\Omega
(\cM)/\cI,Q_{\scriptscriptstyle T\cM}) \stackrel{f^*}{\to}
(C^\infty (\cN), Q_{\scriptscriptstyle \cN})\,,\eeq which induces
a map in cohomology. The conditions of local triviality and gluing by use of
$\exp(\cG')-$valued transition cocycle used in the previous section
gives a natural chain map of complexes
\beq\label{extension_of_based_forms} (\Omega (\cF)_\cG
,Q_{\scriptscriptstyle T\cF}) \to (\Omega
(\cM)/\cI,Q_{\scriptscriptstyle T\cM})\,,\eeq where $\cF$ is the
typical fiber of $p$ (Lemma \ref{gauge_invariance}), the composition
of which with (\ref{I-chain}) determines the characteristic map
(\ref{characteristic_map}) in Theorem
\ref{construction_of_charclasses}.

\vskip 2mm\noindent Let us consider an exact sequence of Lie
algebras \beq\label{exact_sequence_Lie_alg} 0\to
\h\to\g\stackrel{p}{\to}\g_0\to 0\,,\eeq which defines a locally
non-trivial Q-bundle $p\colon \g [1]\to \g_0 [1]$ (we denote the
induced projection by the same letter $p$). The non-triviality
means precisely that the exact sequence
(\ref{exact_sequence_Lie_alg}) does not split in the category of
Lie algebras; here local triviality would imply the global one.
The chain map (\ref{I-chain}) provided by a gauge field $\gf\colon
\g_0\to\g$ in these settings induces a chain map $S^\dot (\h^*)\otimes
\Lambda^\dot (\g^*)\to \Lambda^\dot (\g_0^*)$. Note that here we used
a natural identification induced by the embedding $\iota \colon \h \to \g$.
This chain map, composed
with the embedding $S^\dot (\h^*)^G\to S^\dot (\h^*)\otimes
\Lambda^\dot (\g^*)$, gives nothing but the characteristic
map of exact triples of Lie algebras  of Lecomte \cite{Lecomte}.
For the definition of $S^\dot (\h^*)^G$ one again makes use
of the embedding $\iota$. It is this embedding that induces the
proper replacement of
(\ref{extension_of_based_forms}).

 The complete construction of the
characteristic map of Lecomte  involves twistings by
representations $(V,\rho)$ of $\g_0$, such that the result is taking
values in $H^\dot (\g_0, V)$. In the picture above we need to replace
$\g_0$ and $\g$ with the semidirect products
$\tilde\g_0=\g_0\ltimes_{\rho^*} V^*$ and
$\tilde\g=\g\ltimes_{\rho^*\circ p} V^*$, respectively, which allows
extending the sequence (\ref{exact_sequence_Lie_alg}) canonically to:
\beqn 0\to \h\to\tilde\g\stackrel{\tilde p}{\to}\tilde\g_0\to 0\,.\eeq
Applying the characteristic map $S^\dot (\h^*)^G\to H^\dot
(\tilde\g_0)$ and using the natural isomorphism $H^q
(\tilde\g_0)=\oplus_{k+l=q}H^k (\g_0, \Lambda^l V)$, we immediately
obtain the characteristic classes of Lecomte taking values in all
exterior powers of the representation $\rho$.

\vskip 2mm \noindent It may be interesting to find more general conditions
for a Q-bundle (weaker than those in Definition
\ref{defbundle}) providing an extended version of the map
(\ref{extension_of_based_forms}) which includes simultaneously the
construction of Theorem \ref{characteristic_map} and the Lecomte
characteristic classes.

\section*{Acknowledgments}

\noindent The work of A.K. was partially supported by the research
grant R1F105L10 of the University of Luxembourg. The authors
greatly acknowledge the possibility to work on the completion of
the article at the Erwin Schr\"{o}dinger Institute during the
programm ``Poisson sigma models, Lie algebroids, deformations, and
higher analogues''. T.S.~is grateful to T.~Masson for drawing our
attention to Ref.~[11].



\end{document}